\documentclass{amsart}
\usepackage{amsmath, amsfonts, amssymb, amsthm, amscd}

\newtheorem{theorem}{Theorem}[section]

\theoremstyle{definition}
\newtheorem{definition}[theorem]{{Definition}}
\newtheorem{example}[theorem]{Example}

\newtheorem*{chunk*}{}

\numberwithin{equation}{section}

\theoremstyle{remark}
\newtheorem{remark}[theorem]{Remark}

\newcommand{\hk}{{\rm e}_{HK}} 

\usepackage{graphicx}

\input xy
\xyoption{all}

\newcommand{{\aA}}{\bf A}

\newcommand{\ba}{{\bf a}} 
\newcommand{\bb}{{\bf b}} 
\newcommand{\bt}{{\bf t}}

\newcommand{\sfk}{\mathsf k}

\newcommand{\B}{\mathcal B}

\newcommand{\eI}{\operatorname{extraInequality}}
\newcommand{\iI}{\operatorname{initialInequalities}}

\newcommand{\showtwo}[3]
{%
    \begin{center}
        \includegraphics[width=#3\linewidth]{#1}%
        \hskip1in
        \includegraphics[width=#3\linewidth]{#2}%
    \end{center}
}

\begin{document}
%\subjclass{16Y60; 52B11; 52A07}
\title[Automating the Calculation of Hilbert-Kunz Multiplicities and
F-Signatures]
{Automating the Calculation of Hilbert-Kunz Multiplicities and
F-Signatures}

\author{Gabriel Johnson}
\address{Department of Mathematics, University of Mississippi, Oxford, MS  38677}
\email{gsjohnson@go.olemiss.edu}

\author{Sandra Spiroff}
\address{Department of Mathematics, University of Mississippi, Oxford, MS  38677}
\email{spiroff@olemiss.edu}

\begin{abstract}  The Hilbert-Kunz multiplicity and $F$-signature are important invariants for researchers in commutative algebra and algebraic geometry.  We provide software, and describe the automation of a calculation, for the two invariants
in the case of intersection algebras over polynomial rings.
\end{abstract}

\thanks{
%G.~Johnson's has a degree in Computer Science from the University of Mississippi. S.~Spiroff, Associate Professor of Mathematics, 
S.~Spiroff was supported by Simons Foundation Collaboration Grant \#245926.}

%68W99 is Algorithms

\maketitle

{\small{Code Metadata}}

{\tiny {
\bigskip
  \begin{tabular}{c l l l l l l}
    \hline
    C1 & Current code version & {\text {             }} & {\text {             }}  & v1 \\
    C2 & Permanent link to code/repository & & & https://gitlab.com/gsjohnso/hilbert-kunz-and-f-sig \\
    & used for this code version  \\
    C3 & Legal Code License  & & & BSD 3-Clause \\
    C4 & Code versioning system used & & & Git \\
    C5 & Software code languages, tools, & & & Common Lisp  \\
    & and services used \\
    C6 & Compilation requirements, operating & & & Clozure Common
     Lisp, Bash/Cygwin, Sed \\
    & environments \& dependencies & & & Mathematica (optional); Linux, Mac, Windows  \\
    C7 & If available Link to developer  \\
    & documentation/manual \\
    C8 & Support email for questions & & & gsjohnso@go.olemiss.edu  \\ \hline
  \end{tabular} } }

%%%%%%%%%%%%%%%%%%%%%%%%%%%%%%%%%%%%%%%%%%%%%%%%%

\section{Motivation and Significance}

%%%%%%%%%%%%%%%%%%%%%%%%%%%%%%%%%%%%%%%%%%%%%%%%%

The Hilbert-Kunz multiplicity, along with the $F$-signature, has much importance in the related fields of commutative algebra and algebraic geometry, specifically, in characteristic $p > 0$ methods, but is notoriously difficult to compute in practice. (The article \cite{Hu} by C.~Huneke provides a comprehensive overview of results on these two invariants.) In particular, very few examples exist in the literature where both values are known simultaneously.  This state of affairs motivated the work of F.~Enescu and S.~Spiroff, who calculated the two invariants, as well as the Hilbert-Samuel multiplicity and divisor class group, for certain classes of intersection algebras \cite{ES}.   Because, in their (toric) setting, the Hilbert-Kunz multiplicity and $F$-signature can be realized as volumes of polytopes \cite[Theorem 2.2]{Kz}, \cite[Proposition 4.2]{ES} and \cite[Theorem 3.2.3]{K}, the calculation, being combinatorial in nature, lends itself perfectly to computer automation.  One might expect to be able to obtain general formul{\ae} for the invariants, however, a significant hurdle is the lack of a usable description of the unique Hilbert basis elements in terms of the parameters.  

Let $R$ be a commutative ring and $I, J$ ideals of $R$.  The intersection algebra of $I$ and $J$ is $\mathcal B = \mathcal B_R(I,J) = \bigoplus_{r, s \in \mathbb N}(I^r \cap J^s)$.  In particular, when $R = \sfk[x_1, \dots, x_n]$, for $\sfk$ a field, and $I = (x_1^{a_1}x_2^{a_2}\cdots x_n^{a_n})$ and $J = (x_1^{b_1}x_2^{b_2}\cdots x_n^{b_n})$, with $a_i, b_i \in \mathbb N$, the ordered pairs $(b_i, a_i)$ partition the first quadrant of the plane into a fan made up of pointed rational cones, from which monoids, and Hilbert bases, are obtained.
While for general $n$ and $\ba = (a_1, \dots, a_n), \bb = (b_1, \dots, b_n)$ there is no known formula for $\hk{(\mathcal B)}$ or $s(\mathcal B)$, for any specific numerical entries, Enescu and Spiroff have an algorithm to compute $\hk(\mathcal B)$ and $s(\mathcal B)$.  This algorithm has been automated by G.~Johnson.  See Section 4 for obtaining and using the program.

Intersection algebras overlap with other classes of rings. For example, when $R = \sfk [x]$, $\mathcal B_R((x),(x))$
is a Segre product \cite[\S 3]{W}, i.e., the homogeneous coordinate ring for the Segre embedding $\mathbb P^1 \times \mathbb P^1 \to \mathbb P^3$, and $\mathcal B_R((x^{a_1}),(x))$
is isomorphic to the rational normal scroll $S = \sfk[T, xT, xyT, yT, x^{-1}yT, x^{-2}yT, \dots, x^{-(a_1-1)}yT]$ \cite[Example 3.5]{WY1}.  For details on intersection algebras, especially background material, see \cite[\S 2]{M} and \cite{EM}, and for the complete set of results obtained by Enescu and Spiroff, see \cite{ES}.  References for geometric results are \cite{Hu}, \cite{W}, \cite{WY1}.  Throughout the paper, $R$ will be a polynomial ring over a field and $I, J$ principal monomial ideals.  Since the numerical invariants in our case are given by volumes, they are independent of the characteristic of the field.

%%%%%%%%%%%%%%%%%%%%%%%%%%%%%%%%%%%%%%%%%%%%%%%%%

\section{Description of Software}

%%%%%%%%%%%%%%%%%%%%%%%%%%%%%%%%%%%%%%%%%%%%%%%%%

We extract the necessarily positive\footnote{The case of non-negative integers can be addressed within the scope of positive integers.  See \cite[Proposition 1.6]{ES}.} exponents into two strings of positive integers, $\ba = a_1, \dots, a_n$ and $\bb = b_1, \dots, b_n$ since they may be treated without any reference to the ideals themselves.  The command line interface can be run on a Unix shell, i.e. on Linux, Mac, or on Windows using Cygwin. For the parameters, the user need only enter the two
lists of exponents delimited by spaces one after the other, with no delimiter
marking the end of one list and the beginning of the next; the list of numbers
entered is automatically cut in half and the first half assumed to be the
exponents of one ideal and the second half the exponents of the other.

\begin{tabular}{l l l r}
{\it user directory}\$ {\tt ./calculate-integral} $a_1$ $a_2 \cdots a_n$ $b_1$ $b_2 \cdots b_n$ \\
\end{tabular}

The program displays the values of the Hilbert-Kunz multiplicity and the $F$-signature on two lines and terminates, as shown below.

\begin{tabular}{l l l r}
 {\tt Hilbert-Kunz Multiplicity} & & & $= p_1/q_1$ \\
 {\tt F-Signature} & & & $= p_2/q_2$  \\
\end{tabular}

When the {\tt ./inequalities} command is used, as shown below, the program displays the {\it Mathematica} code used to calculate the integrals.

\begin{tabular}{l l l r}
{\it user directory}\$ {\tt ./inequalities} $a_1$ $a_2 \cdots a_n$ $b_1$ $b_2 \cdots b_n$ \\
\end{tabular}

\subsection{Software architecture and functionality}

Three main technologies are used.

1. The Bash shell. This collects the arguments from the user and
    invokes the program itself.

2.  Clozure Common Lisp. This is the programming language we used to
    calculate the set of inequalities bounding the solid in question. It
    provides a large number of list-processing operators as well as highly
    flexible looping constructs, making it well suited for dealing with the
    points of indefinite dimension involved in this project.
 
 3. {\it Mathematica} \cite{Wolf}. This is used to calculate the final values of the
    integrals that give the Hilbert-Kunz multiplicity and $F$-signature, in any
    number of dimensions.  The integration features of {\it Mathematica} are
    particularly well suited to this project, where the region of integration is
    defined by a set of inequalities in any number of dimensions.
    Other tools capable of solving integrals tend to require a function rather
    than a set of inequalities to define the area of integration, and limit the
    number of dimensions to two or three.

The terminal interface consists of a directory, called ccl, containing Clozure
Common Lisp itself, a directory, called src, containing the Lisp source code as
well as the Sed script used to format the output readably, and two bash scripts,
setup and calculate-integral. The first, setup, must be run before attempting to
run the other. It has Clozure compile the source code into an image file
(effectively a DLL; in Lisp, all programs are essentially implemented as a DLL
loaded by the language kernel with a specified top-level function to run instead
of a REPL) called hkm.image. Once hkm.image exists in the same directory as the
other files, calculate-integral can be run as specified above.

All the files and directories must occupy the same directory in
order for the program to work. Further, the program also assumes that the
Mathematica kernel, MathKernel, exists in a specified location depending on the
operating system. 

\subsection{Implementation details}  Below is a brief overview of the algorithm used in our program.  Concurrently, we include a specific example to demonstrate the mathematics at each (hidden) step.

Recall that, to implement the program, $a_1, a_2, \dots, a_n$ is entered as the exponents of one of the ideals, followed by $b_1, b_2, \dots, b_n$ for the other, {\it with the lists, given or necessarily permuted, in fan order}.  Once the arguments are collected from the command line, all steps but the last are done entirely in Lisp.

\begin{example}
For $R = \sfk[x,y]$ and $I = (x^5y^2), J = (x^2y^3)$, enter the command {\tt ./calculate-integral 5 2 2 3} in the Bash shell.
\end{example}

{\bf Step 1.}  Our program first finds the Hilbert set.  As in \cite[\S 1]{ES}, \cite{M}, there are cones $C_i$ and $C$, and monoids $Q_i$ associated to the pairs $(b_i, a_i)$, and the unique Hilbert basis $\mathcal H_i$ for each $Q_i$ may be found via the algorithm given in \cite[Algorithm 2.4]{CDE}.  

\begin{definition}\cite[Definition 1.2]{ES} \label{Hilbertsets} The set
  $\mathcal H= \cup _{i=0, \dots, n} \mathcal H_i$ is called the {\bf Hilbert
  set} for $\B(I, J)$.  For $v=(r,s) \in \mathcal H$, 
let $\bt(v) = (\max(a_ir, b_is))_{i=1, \ldots, n}$.
  Set $\mathcal G= \{ (v, \bt(v)) : v \in \mathcal H\}$.
\end{definition}

e.g., $\mathcal H_0 = \{(0,1), (1,3), (2,5)\},$ $\mathcal H_1 = \{(1,1), (1,2), (2,5), (3,2)\}$, and 
$\mathcal H_2 = \{(1,0),$ \\
$ (2,1), (3,2)\},$ hence $|\mathcal H| = |\mathcal H_0 \cup \mathcal H_1 \cup \mathcal H_2| = 8$.

{\bf Step 2.} It finds the set $\mathcal G$. 
Each pair in the Hilbert set will be rewritten as an $(n+2)$-tuple.
The extra $n$ coordinates $z_i$ are determined  according to the formula $z_i(x,
y) = \max(a_ix, b_iy)$, corresponding to the intersection of monomial ideals.

e.g., $\mathcal G = \{(1,0,5,2), (0,1,2,3), (2,5,10,15), (3,2,15,6), (2,1,10,4), (1,2,5,6),$\\
$ (1,1,5,3), (1,3,6,9) \}$

{\bf Step 3.} It determines the set of
inequalities bounding a solid in $n+2$ dimensions. If the individual
coordinates of each element of $\mathcal G$ are labelled as $(x, y, z_1, z_2, \dots,
z_n)$, then an initial set of inequalities is derived by a predetermined
template: 

\centerline{$\displaystyle{\bigwedge_{i=1}^n 0 \leq a_ix \leq z_i \land 0 \leq b_iy \leq z_i.}$}

From these initial inequalities and $\mathcal G$, further
inequalities must be derived, and the conjunction of all these inequalities will produce the polytope we seek: an element $(p, q, r_1, r_2, \dots, r_n) \in \mathcal G$ generates a set of further
inequalities by the following function: \hskip.25in $\eI(p, q, r_1, r_2, \dots, r_n) = $

\centerline{$\neg \iI(x - p, y - q, z_1 - r_{1}, z_2 - r_{2}, \dots, z_{n} - r_{n}),$}
where $\iI$ is a function that accepts an element of $\mathcal G$ and returns true if
and only if the given point satisfies all the initial inequalities.  The program repeats this process for every element of $\mathcal G$.
%This process is repeated for every element of $\mathcal G$.
%If $\iI$ is a function that accepts a point and returns true if and only if the
%point satisfies the initial inequalities and $\eI$ is a function that accepts a
%point and returns true if and only if the point satisfies the inequality to be
%newly derived, then $\eI$ can be defined in terms of $\iI$ as $$\eI(p, q, r_1,
%r_2, \dots, r_n) = $$ $$\neg \iI(x - p, y - q, z_1 - r_{1}, z_2 - r_{2}, \dots,
%z_{n} - r_{n}).$$
Then, it joins all the sets of inequalities by conjunction, producing our polytope.  Its volume is the Hilbert-Kunz multiplicity \cite[Theorem 2.2]{Kz}, \cite[Proposition 4.2]{ES}, and will appear in the first line of output.

e.g., {\small{  {\tt Integrate[Boole[((x >= 0 \&\& y >= 0 \&\& z >= 5*x \&\& z >= 2*y \&\& } \\
{\tt w >= 2*x \&\& w >= 3*y)) } \\
{\tt   \&\& ((z < 5*x + 1 || z < 2*y + 1)) } \\
{\tt    \&\& ((w < 2*x + 1 || w < 3*y + 1)) } \\
{\tt    \&\& ((x < 1 || z < 2*y + 5 || w < 3*y + 2)) } \\
{\tt    \&\& ((y < 1 || z < 5*x + 2 || w < 2*x + 3)) } \\
{\tt    \&\& ((x < 2 || y < 5 || w < 2*x + 11)) } \\
{\tt    \&\& ((x < 3 || y < 2 || z < 2*y + 11)) } \\
{\tt    \&\& ((x < 2 || y < 1 || z < 2*y + 8 || w < 3*y + 1)) } \\
{\tt    \&\& ((x < 1 || y < 2 || z < 2*y + 1 || w < 2*x + 4)) } \\
{\tt    \&\& ((x < 1 || y < 1 || z < 2*y + 3 || w < 2*x + 1)) } \\
{\tt    \&\& ((x < 1 || y < 3 || z < 5*x + 1 || w < 2*x + 7))], } \\
 {\tt  \{x, 0, 100\}, \{y, 0, 100\}, \{z, 0, 500\}, \{w, 0, 500\}] \} } }}

{\bf Step 4.} It determines the inequalities for the $F$-signature, which is
much simpler, by definition, as they follow the formula 

\centerline{$\displaystyle{\bigwedge_{i = 1}^n
a_ix \leq z_i \leq 1 + a_ix \land b_iy \leq z_i \leq 1 + b_iy.}$} 
This can be found with a simple for loop.  The $F$-signature is the volume of this second
polytope \cite[Theorem 3.2.3]{K}, and will appear in the second line of output.

e.g.,{\small{  {\tt Integrate[Boole[5*x <= z < 1 + 5*x \&\& 2*y <= z < 1 + 2*y \&\& }
  {\tt 2*x <= w < 1 + 2*x \&\& 3*y <= w < 1 + 3*y ], \{x, 0, 1\}, \{y, 0, } 
{\tt  1\}, \{z, 0, 25\}, \{w, 0, 15\}] } }}

If the {\tt ./inequalities} command was used, then the program terminates after displaying the inequalities.  Otherwise, the program passes the complete sets of inequalities to {\it Mathematica}, which calculates the volumes of the corresponding solids.
Sed is used to format the output readably.

\begin{tabular}{l l l r}
{\tt Hilbert-Kunz Multiplicity} & & & {\tt = 37283/9900} \\
{\tt F-Signature} & & & {\tt = 1087/29700} \\
\end{tabular}

Implementation of the {\tt ./inequalities} command allows the user to bypass the use of {\it Mathematica}. This allows one to obtain 
the inequalities in Steps 3 and 4.  Not only is this useful for 3D graphing and printing, but the user may integrate or analyze via another program, especially if he/she does not have a local installation of {\it Mathematica}, a proprietary software, or only wants to study the $F$-signature.  As evinced in the examples below, the running time in higher dimensions can significantly increase depending upon the complexity of the
geometry involved.  However, the bottleneck in these cases is the integration related to the Hilbert-Kunz multiplicity; the steps calculating the regions determined by the inequalities, and the integration related to the $F$-signature, finished almost instantaneously even in the most complex cases tested.

%%%%%%%%%%%%%%%%%%%%%%%%%%%%%%%%%%%%%%%%%%%%%%%%%

\section{Some Examples-Polytopes and Running Times}

%%%%%%%%%%%%%%%%%%%%%%%%%%%%%%%%%%%%%%%%%%%%%%%%%

\begin{example} (Volumes.) The 3D polytopes
associated to the invariants for ${\mathcal B}_R((x^3), (x^2))$, where $R = \sfk[x]$, show the volumes calculated by the program.
   \begin{figure}[tbh] \showtwo{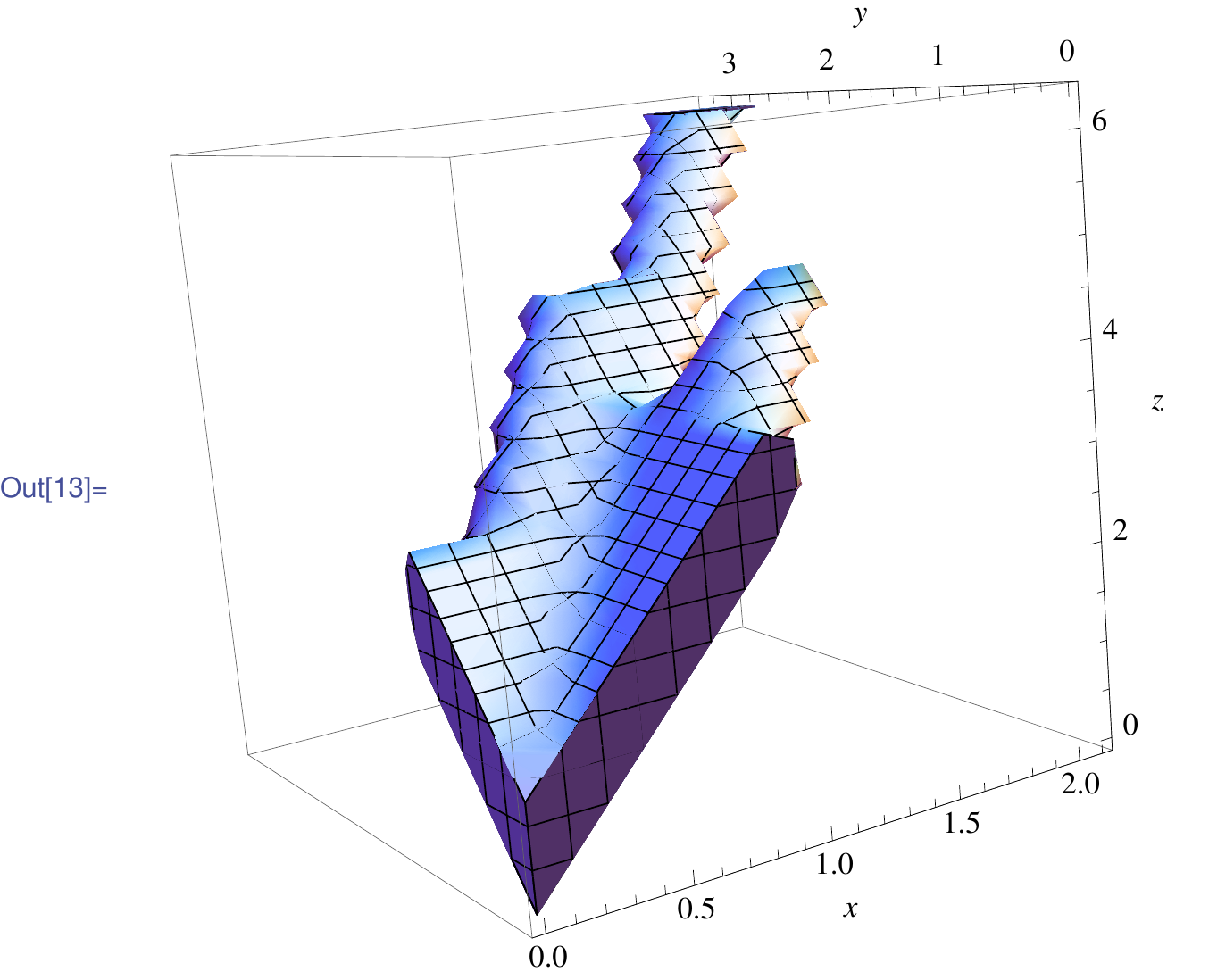}{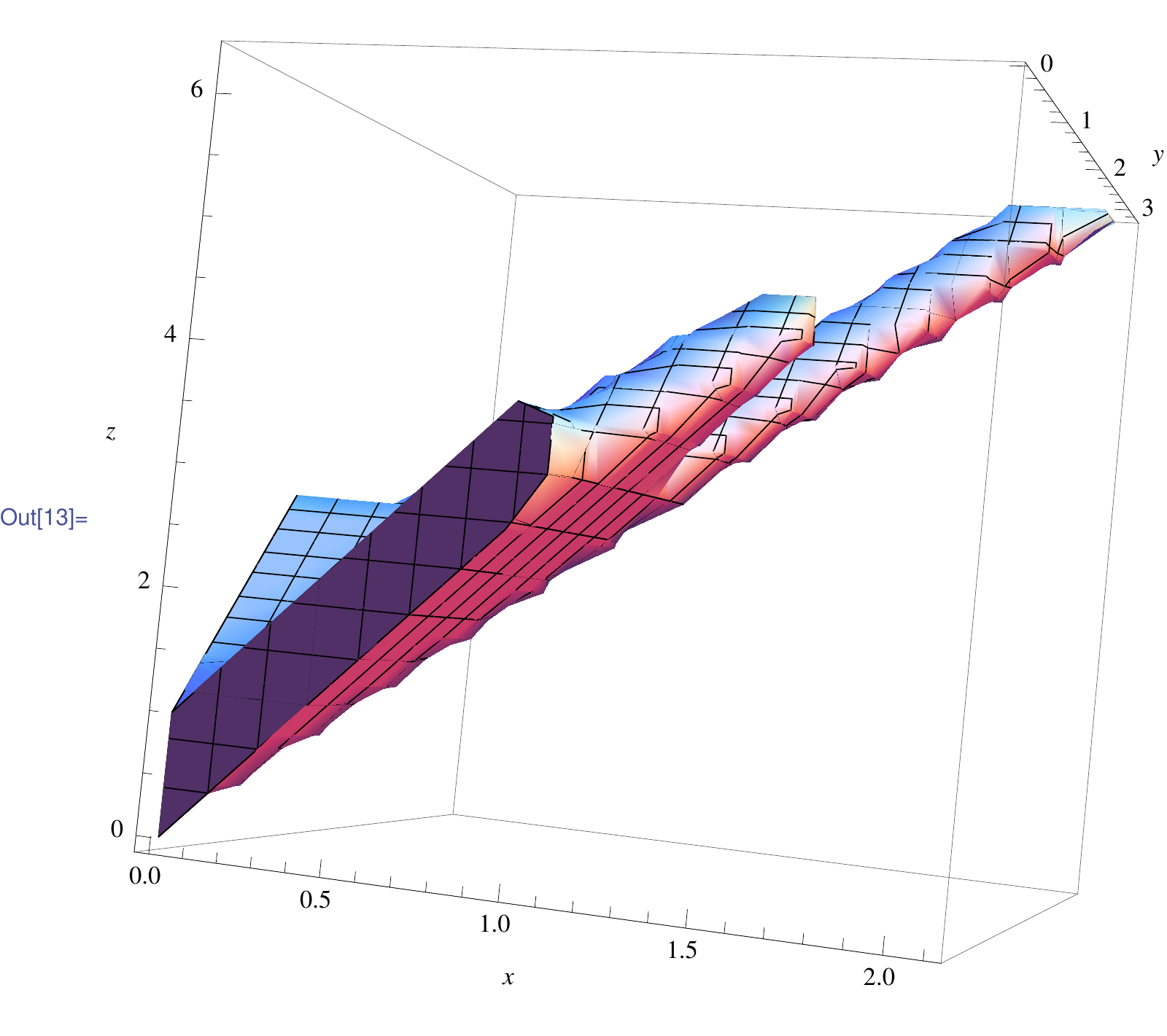}{0.25} \caption{{\bf Hilbert
  Kunz Multiplicity}: view from the front and from the side, $\hk({\mathcal B}) = \frac {41}
  {18}.$} \end{figure}
  \begin{figure}[tbh] \showtwo{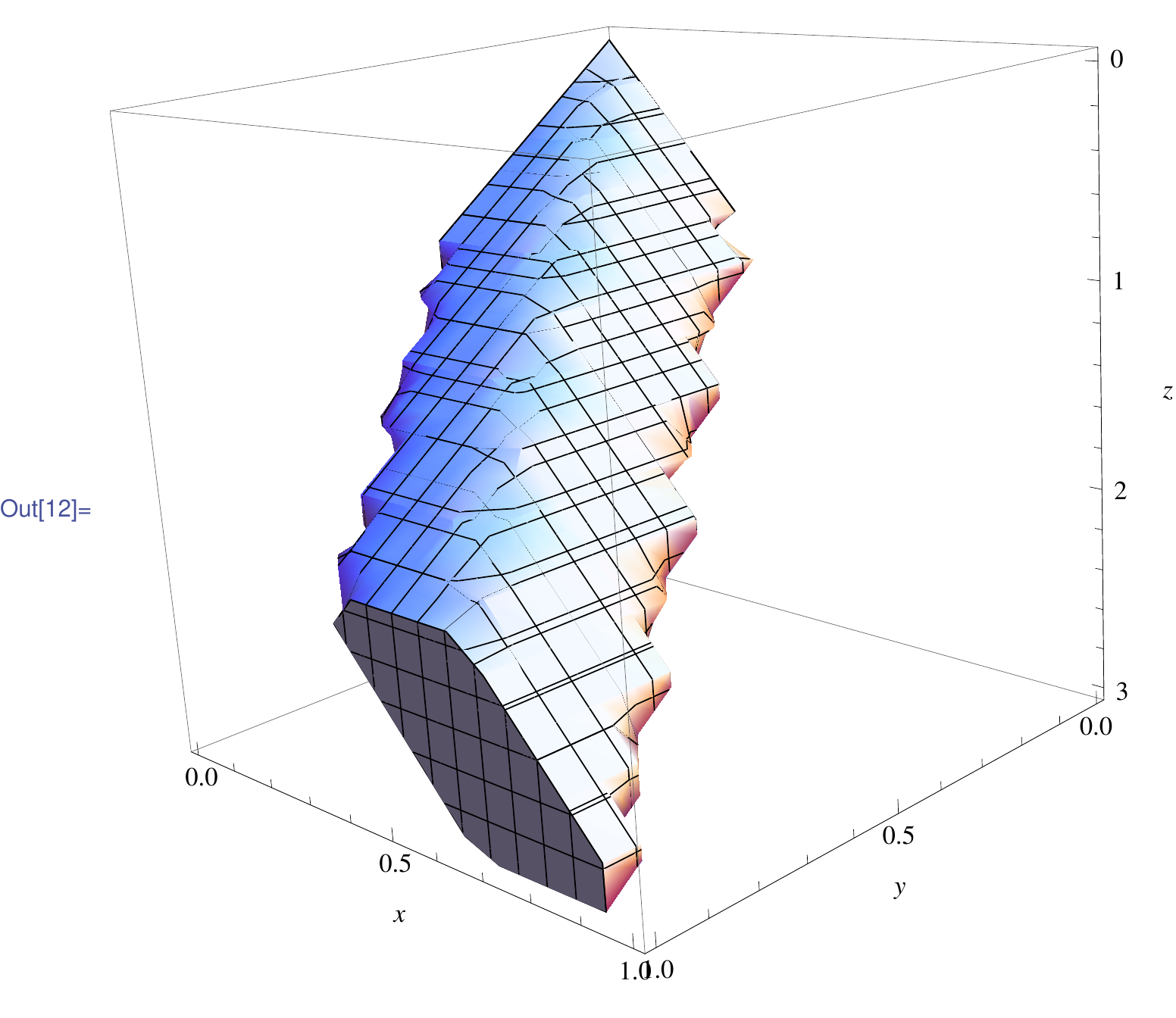}{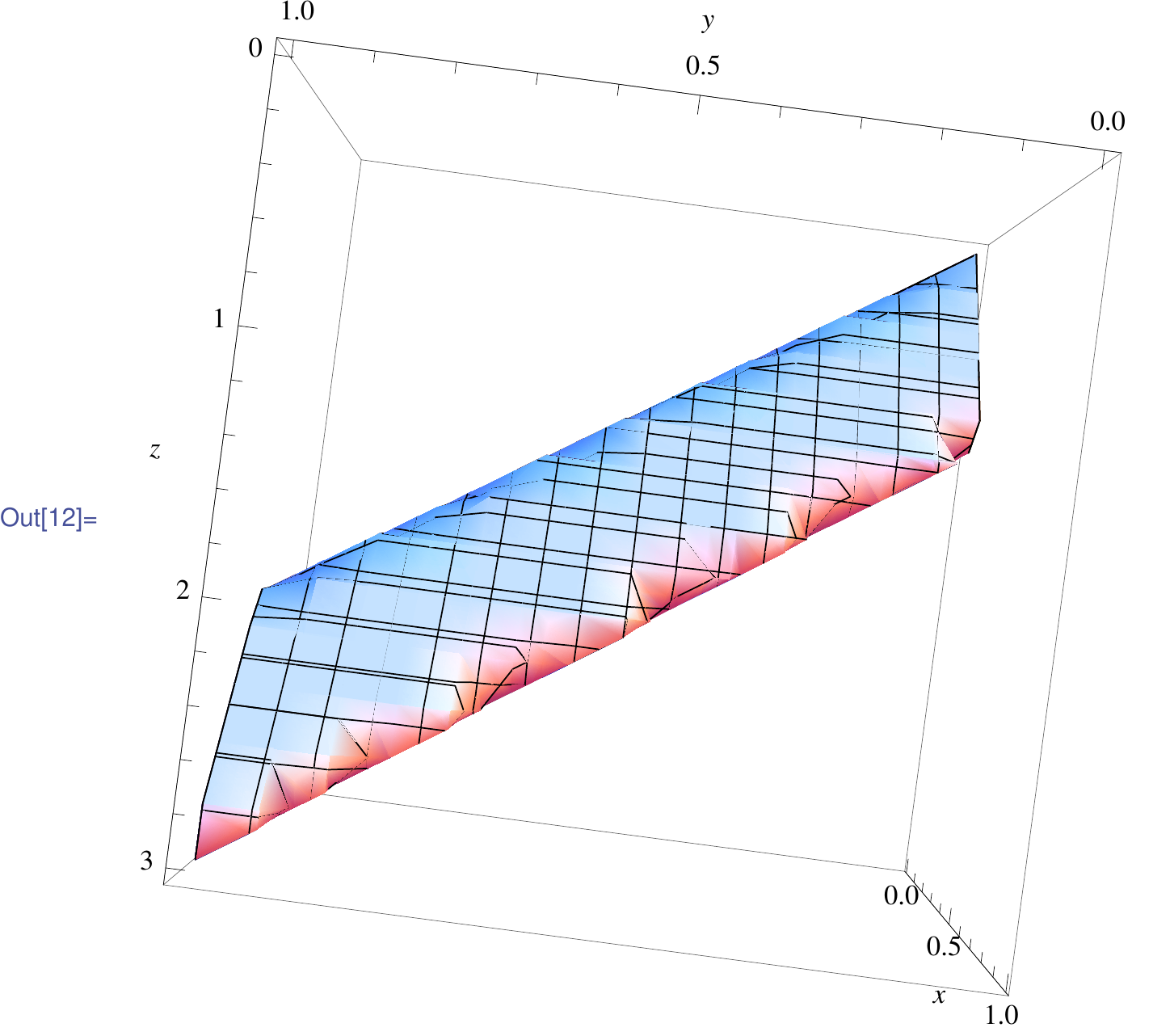}{0.25} \caption{{\bf
  $F$-Signature}: view from the front and from the side, $s(\mathcal B) = \frac {11} {36}.$} \end{figure}
\end{example}

\begin{example} (Fan order and Running times.) \label{5D} When $R = \sfk[x,y,z]$ and $I = (x^2y^4z^7)$, $J = (xy^6z^5)$, the calculation had a running time of approximately 3.5 hours on a laptop, although the inequalities command outputs almost immediately.

\begin{tabular}{l l l r}
{\tt ./calculate-integral 2 7 4 1 5 6} \\
{\tt Hilbert-Kunz Multiplicity} & & & {\tt = 1874881259711/391184640000} \\
{\tt F-Signature} & & & {\tt = 27251293/1564738560} \\
\end{tabular}

On the other hand, the example {\tt ./calculate-integral 1 6 5 1 6 5} executes in a matter of seconds; i.e., $\hk({\mathcal B}) = \frac {1633}{864}, s({\mathcal B}) = \frac {95}{864}$, for $\mathcal B = \mathcal B_R(J, J)$.
%\vskip.2in

This last calculation illustrates the fact that if $a_i = b_i$ for all $i$, then $\hk{(\mathcal B)} + s(\mathcal B) = 2$.  (See \cite[Proposition 4.3]{ES}.)
\end{example}

%%%%%%%%%%%%%%%%%%%%%%%%%%%%%%%%%%%%%%%%%%%%%%%%%

\section{Software Availability and Directions for Use}

%%%%%%%%%%%%%%%%%%%%%%%%%%%%%%%%%%%%%%%%%%%%%%%%%

The software may be obtained at either of the websites: \begin{verbatim}http://home.olemiss.edu/~spiroff/ \end{verbatim} or 

\begin{verbatim}https://gitlab.com/gsjohnso/hilbert-kunz-and-f-sig \end{verbatim}
\bigskip

\underline{One time steps}

1. Download the .zip file and extract it into your chosen directory. 

2. In a Bash shell, go to the chosen directory using the {\tt cd} command(s) and run {\tt ./setup}. 
\smallskip

\underline{Commands for regular use}-in a Bash shell, in the chosen directory (in fan order):

 {\tt ./calculate-integral}  $a_1$ $a_2 \cdots a_n$ $b_1$ $b_2 \cdots b_n$

 {\tt ./inequalities}  $a_1$ $a_2 \cdots a_n$ $b_1$ $b_2 \cdots b_n$
\bigskip

\begin{remark} The details behind this automation, and the mathematical theory regarding the computations, are in the paper {\it Computing the invariants of intersection algebras of principal monomial ideals}, by F.~Enescu and S.~Spiroff, International Journal of Algebra and Computation, to appear.
\end{remark}

\section*{Acknowledgements}
The authors thank Florian Enescu at Georgia State University, Dawn Wilkins and Conrad Cunningham in the Department of Computer Science at the University of Mississippi, and the various specialists in the Mississippi Center for Supercomputing Research on the Oxford campus for useful conversations about this project.

\end{document}